\documentclass[12pt,cs4size,a4paper,fancyhdr,fntef]{article}
\ifx\pdfoutput\undefined 
  \usepackage[dvipdfm,CJKbookmarks]{hyperref}
\else
  \usepackage[pdftex,CJKbookmarks]{hyperref}
\fi
\usepackage{shortvrb,ulem,makeidx} 
\usepackage{fancyhdr}
\usepackage{multicol}
\usepackage{array,color}
\usepackage{amssymb,amsthm}
\usepackage{amsfonts}
\usepackage{extarrows}
\usepackage{graphicx}
\usepackage{makecell,rotating}
\usepackage{multirow}
\usepackage{url}
\usepackage{indentfirst}
\usepackage{subfigure}
\usepackage{mathrsfs}

 \MakeShortVerb{\|}

 \numberwithin{equation}{section}

\newtheorem{defi}{DEFINITION}[section]
\newtheorem{theo}[defi]{THEOREM}
\newtheorem{prop}[defi]{PROPOSITION}
\newtheorem{lemm}[defi]{LEMMA}

\begin{document}

\centerline{\large\bf Conditional G-expectation in $\mathbb{L}^{p}$ and Related It\^o's Calculus }

\vskip 0.6cm

{\bf Y. FAN\footnote[0]{This research was supported by Beijing Natural Science Foundation Grant No. 1112009.\\
\indent\ \  Address correspondence to Y. Fan, School of Science, North China University of Technology, Beijing 100144, China; Email: fanyl@ncut.edu.cn} }

\vskip0.2cm

School of Science, North China University of Technology,   China

\date{}

\vskip 0.3cm

In this paper, we define a dynamically consistent conditional G-expectation in space $\mathbb{L}^{p}$,
and give the related stochastic calculus of It\^o's type, especially get It\^o's formula for a general $C^{1,2}$-function.

\vskip 0.3cm


\vskip 0.3cm


\vskip 0.3cm

\section{Introduction}

The notion of G-expectation, which
is a typical nonlinear expectation, was proposed by Peng \cite{Peng (2006)}, \cite{Peng (2008)}, and \cite{Peng (2010)}. It can be regarded as a nonlinear
generalization of Wiener probability space $(\Omega, \mathcal{F}, P)$,
where $\Omega=C([0,\infty),\mathbb{R}^{d})$, $\mathcal{F}=\mathcal{B}(\Omega)$,
and $P$ is a Wiener probability measure defined on $(\Omega, \mathcal{F})$. On the same canonical space $\Omega$, G-expectation is a
sublinear expectation, such that the same canonical process $B_{t}(\omega):=\omega_{t}, t\geq 0$
is a G-Brownian motion, i.e., it is a continuous process with stable and independent
increments. One important feature of G-expectation is its time consistency. To be precise,
let $\xi$ be a random variable and $Y_{t}:=\mathbb{E}^{G}_{t}[\xi]$ be the conditional G-expectation , then one
has $\mathbb{E}^{G}_{s}[\xi]=\mathbb{E}^{G}_{s}[\mathbb{E}^{G}_{t}[\xi]]$
 for any $s<t$. For this reason, the conditional G-expectation is called
a G-martingale, or a martingale under G-expectation.

A well-known and fundamentally important fact in probability theory is that the
linear space $L^{1}_{P}$ coincides with the $E_{P}[\vert\cdot\vert]$-norm completion of the space of
bounded and $\mathcal{F}$-measurable functions
$B_{b}(\Omega)$, or bounded and continuous functions $C_{b}(\Omega)$, or even smaller one, the space $Lip_{b.cyl}(\Omega)\subset C_{b}(\Omega)$ of bounded and Lipschitz
cylinder functions. While in the theory of G-expectation, Denis, et al.\cite{Denis et al. (2011)} proved that
the $\mathbb{E}[\vert\cdot\vert]$-norm completion of $Lip_{b.cyl}(\Omega)$ and $C_{b}(\Omega)$
are the same space $L^{1}_{G}(\Omega)$, the random variables $X=X(\omega)$ which are
quasi-continuous with respect to the natural Choquet capacity $c(A):=\mathbb{E}^{G}[\mathbf{I}_{A}]$, $A\in \mathcal{B}(\Omega)$,
but they are strict subspace of the $\mathbb{E}[\vert\cdot\vert]$-norm completion of $B_{b}(\Omega)$.
Moreover the latter one is, again, a strict subspace of $\mathbb{L}^{1}$, the space of all $\mathcal{F}$-measurable
random variables $X$ such that $\mathbb{E}[\vert X\vert]<\infty$.

In Peng \cite{Peng (2006)}, \cite{Peng (2008)}, and \cite{Peng (2010)}, G-expectation and the related
It\^o's calculus are mainly based on space $L^{1}_{G}(\Omega)$. Then Li and Peng \cite{Li  Peng (2009)} extend the It\^o's integral
to space without the quasi-continuity, obtain It\^o's
integral on stopping time interval, and get It\^o's formula for a general $C^{1,2}$-function, which generalizes the
previous results of Peng \cite{Peng (2006)}, \cite{Peng (2008)}, and \cite{Peng (2010)} and it’s improved version of Gao \cite{Gao (2009)} and Zhang et al. \cite{Zhang et al. (2010)}.

But in Li and Peng \cite{Li  Peng (2009)}, the conditional G-expectation has not been defined,
so whether the martingale properties still hold for the stochastic integral $(\int_{0}^{t}\eta_{s}dB_{s})_{0\leq t\leq T}$ and
the conditional G-expectation of random variables without quasi-continuous condition is open.

This paper is organized as follows. In section 2, we give the definition of conditional G-expectation of random variables in space $\mathbb{L}^{1}$.
In section 3, we define the related It\^o's integral in space $\mathbb{M}^{2}(0,T)$. In section 4, we prove the It\^o's formula for general $C^{1,2}$ function.

\section{Conditional G-expectation in $\mathbb{L}^{p}$}

\subsection{G-Brownian motion and G-expectation}

We first present some preliminaries in the theory of G-expectation and the related space of random
variables. More relevant details can be found in Peng \cite{Peng (2006)}, \cite{Peng (2008)}, \cite{Peng (2010)}.

Let $\Omega$ be a given set and let $\mathcal{H}$ be a linear space of
real functions defined on $\Omega$ such that $c\in\mathcal{H}$ for
each constant c and $\vert X\vert\in\mathcal{H}$ if
$X\in\mathcal{H}$. A sublinear expectation $\mathbb{E}$ on
$\mathcal{H}$ is a functional $\mathbb{E}: \mathcal{H}\rightarrow
\mathbb{R}$ satisfying the following properties: for all $X,Y \in
\mathcal{H}$,

(a) Monotonicity: If $X\geq Y$ then
$\mathbb{E}[X]\geq\mathbb{E}[Y]$.

(b) Constant preserving: $\mathbb{E} [c] = c$.

(c) Sub-additivity: $\mathbb{E}[X+Y]\leq\mathbb{E}[X]+\mathbb{E}[Y
]$.

(d) Positive homogeneity: $\mathbb{E}[\lambda X]
=\lambda\mathbb{E}[X], \forall \lambda\geq 0 $.

The triple $(\Omega,\mathcal{H},\mathbb{E})$ is called a sublinear
expectation space. In the literature, for technical convenience,
$\mathcal{H}$ is taken as the space satisfying that if
$X_{1},\ldots, X_{n}\in \mathcal{H}$ then $\varphi(X_{1},\ldots,
X_{n})\in \mathcal{H}$ for each $\varphi \in
C_{l,Lip}(\mathbb{R}^{n})$ where $C_{l,Lip}(\mathbb{R}^{n})$ denotes
the linear space of (local Lipschitz) functions $\varphi$ satisfying
$$
\mid\varphi(x)-\varphi(y)\mid\leq C(1+\vert x\vert^{m}+\vert y\vert^{m})\mid x-y\mid,
\forall x,y \in \mathbb{R}^{n},
$$
for some $C>0, m \in \mathbb{N}$ depending on $\varphi$. The linear
space $C_{l,Lip}(\mathbb{R}^{n})$ can be replaced by
$\mathbb{L}^{\infty}(\mathbb{R}^{n}), C_{b}(\mathbb{R}^{n}),
C^{k}_{b}(\mathbb{R}^{n}),  C_{unif}(\mathbb{R}^{n}),
C_{b,Lip}(\mathbb{R}^{n})$, and $L^{0}(\mathbb{R}^{n})$.
In this case $X=(X_{1},\ldots,X_{n})$ is called an $n$-dimensional random vector, denoted by $X\in\mathcal{H}^{n}$.

\begin{defi}
In a nonlinear expectation space $(\Omega, \mathcal{H},\mathbb{E})$, a random vector
$Y\in \mathcal{H}^{n}$ is said to be independent from another random vector $X\in \mathcal{H}^{n}$ under $\mathbb{E}[\cdot]$ if
for each test function $\varphi\in C_{l.L_{ip}}(\mathbb{R}^{m+n})$ we have
$$
\mathbb{E}[\varphi(X,Y )]=\mathbb{E}[\mathbb{E}[\varphi(x,Y )]_{x=X}].
$$
\end{defi}

Let $\Omega=C^{d}_{0}(\mathbb{R}^{+})$ be the space of all
$\mathbb{R}^{d}$-valued continuous paths $(\omega_{t})_{t\in
\mathbb{R}^{+}}$, with $\omega_{0}=0$, equipped with the distance
$$
\rho(\omega^{1},\omega^{2}):=\sum_{i=1}^{\infty}2^{-i}[(\max_{t\in[0,i]}\vert\omega^{1}_{t}-\omega^{2}_{t}\vert)\wedge
1].
$$
For each fixed $T\in[0,\infty)$, we set
$\Omega_{T}:=\{\omega_{\cdot\wedge T}: \omega\in\Omega\}$.

Let
$$
L_{ip}(\Omega_{T}):=\{\phi(B_{t_{1}\wedge T},\cdots,B_{t_{n}\wedge
T}): n\in\mathbb{N},t_{1},\cdots,t_{n}\in[0,\infty),\phi\in
C_{l,Lip}(\mathbb{R}^{d\times n})\}
$$
with $B_{t}=\omega_{t}, t\in[0,\infty)$ for $\omega\in\Omega$, and
$$
L_{ip}(\Omega):=\bigcup_{n=1}^{\infty} L_{ip}(\Omega_{n}).
$$

Let $(\xi_{i})_{i=1}^{\infty}$ be a sequence of d-dimensional random vectors
on a sublinear expectation space
$(\tilde{\Omega}, \tilde{\mathcal{H}},\tilde{E})$ such that
$\xi_{i}$ is G-normal distributed and $\xi_{i+1}$ is independent
from $(\xi_{1}, \cdots, \xi_{i})$ for each $i=1,2,\cdots$. For each
$X\in L_{ip}(\Omega)$ with
$$
X=\phi(B_{t_{1}}-B_{t_{0}},B_{t_{2}}-B_{t_{1}},\cdots,B_{t_{n}}-B_{t_{n-1}}),
$$
some $\phi\in C_{l,Lip}(R^{d\times n})$ and
$0=t_{0}<t_{1}<\cdots<t_{n}<\infty$, define G-expectation $\hat{\mathbb{E}}[\cdot]$ as
$$
\hat{\mathbb{E}}[\phi(B_{t_{1}}-B_{t_{0}},B_{t_{2}}-B_{t_{1}},\cdots,B_{t_{n}}-B_{t_{n-1}})]:=\tilde{\mathbb{E}}[\phi(\sqrt{t_{1}-t_{0}}\xi_{1},\cdots,\sqrt{t_{n}-t_{n-1}}\xi_{n})].
$$

And the related conditional G-expectation of
$$
X=\phi(B_{t_{1}},B_{t_{2}}-B_{t_{1}},\cdots,B_{t_{n}}-B_{t_{n-1}})
$$
under $\Omega_{t_{j}}$ is defined by
\begin{eqnarray*}
\begin{split}
\hat{\mathbb{E}}[X\vert\Omega_{t_{j}}]:=\psi(B_{t_{1}},\cdots,B_{t_{j}}-B_{t_{j-1}}),
\end{split}
\end{eqnarray*}
where
$$
\psi(x_{1},\cdots,x_{j})=\tilde{\mathbb{E}}[\phi(x_{1},\cdots,x_{j},\sqrt{t_{j+1}-t_{j}}\xi_{j+1},\cdots,\sqrt{t_{n}-t_{n-1}}\xi_{n})].
$$
$\hat{\mathbb{E}}[\cdot]$ consistently defines a sublinear
expectation on $L_{ip}(\Omega)$ and $(B_{t})_{t\geq0}$ is a
G-Brownian motion.

The sublinear expectation $\hat{\mathbb{E}}[\cdot]: L_{ip}(\Omega)\rightarrow \mathbb{R}$ defined through the above procedure is called a
G-expectation. The corresponding canonical process $(B_{t})_{t\geq 0}$ on the sublinear expectation space $(\Omega, L_{ip}(\Omega),\hat{\mathbb{E}})$
is called a G-Brownian motion.

Let $L^{p}_{G}(\Omega_{T}), p\geq 1$, denotes the completion of
$$
L_{ip}(\Omega_{T}):=\{\varphi(B_{t_{1}\wedge
T},\cdots,B_{t_{n}\wedge T}): n\in
\mathbb{N},t_{1},\cdots,t_{n}\in[0,\infty),\varphi\in
C_{l,Lip}(\mathbb{R}^{d\times n})\}
$$
under the norm $\parallel X\parallel_{p}:=(\hat{\mathbb{E}}[\vert
X\vert^{p}])^{\frac{1}{p}}$. And set $L_{ip}(\Omega):=\bigcup_{n=1}^{\infty}L_{ip}(\Omega_{n})$.
The definition of G-expectation and conditional G-expectation can be extended to space $L^{p}_{G}(\Omega)$.

\subsection{Conditional G-expectation in $\mathbb{L}^{p}$}

Let $\mathcal{M}$ be the collection of all probability measures on $(\Omega, \mathcal{B}(\Omega))$, and

$L^{0}(\Omega)$: the space of all $\mathcal{B}(\Omega)$-measurable real functions;

$L^{0}(\Omega_{t})$: the space of all $\mathcal{B}(\Omega_{t})$-measurable real functions;

$B_{b}(\Omega)$: all bounded elements in $L^{0}(\Omega)$; $B_{b}(\Omega_{t})$: all bounded elements in $L^{0}(\Omega_{t})$;

$C_{b}(\Omega)$: all bounded and continuous elements in $L^{0}(\Omega)$; $C_{b}(\Omega_{t})$: all bounded and continuous elements in $L^{0}(\Omega_{t})$.

\begin{theo}(Denis, et al.\cite{Denis et al. (2011)})\label{representation theorem of G-expectation} There exists a weakly compact
subset $\mathscr{P}\subseteq \mathcal{M}$, such that
\begin{eqnarray*}
\hat{\mathbb{E}}[\xi]=\sup_{P\in\mathscr{P}}E_{P}[\xi], \forall \xi\in L^{1}_{G}(\Omega).
\end{eqnarray*}
$\mathscr{P}$ is called a set that represents $\hat{\mathbb{E}}$.
\end{theo}

The upper expectation
of probability measure set $\mathscr{P}$ is defined in Huber and Strassen \cite{Huber (1973)}:
For
each $X\in L^{0}(\Omega)$ such that $E_{P}[X]$ exists for each
$P\in\mathscr{P}$, the upper expectation about  $\mathscr{P}$ is
defined as
$$
\mathbb{E}[X]=\mathbb{E}^{\mathscr{P}}[X]:=\sup_{P\in\mathscr{P}}E_{P}[X].
$$
Denote
$$
c(A)=\sup_{P\in\mathscr{P}}P(A), A\in \mathcal{B}(\Omega).
$$
Then $c(\cdot)$ is a Choquet capacity. A set A is called polar if
$c(A)=0$, and a property holds ``quasi-surely''(q.s.) if it holds
outside a polar set. Let
\begin{eqnarray*}
\begin{split}
&\mathcal{L}^{p}:=\{X\in L^{0}(\Omega): \mathbb{E}[\vert X\vert^{p}]<\infty\}, 0<p<\infty,\\
&\mathcal{L}^{\infty}:=\{X\in L^{0}(\Omega): \exists \text{ a constant } M, s.t. \vert X\vert\leq M, q.s.\}, \\
&\mathcal{N}:=\{X\in L^{0}(\Omega):  X=0, c-q.s.\},\\
&\mathbb{L}^{p}:=\mathcal{L}^{p}/\mathcal{N}.
\end{split}
\end{eqnarray*}
Then for $0<p<1$, $\mathbb{L}^{p}$ is a complete metric space under
the distance $d(X,Y):=\mathbb{E}[\mid X-Y\mid^{p}]$,
and for $1\leq p<\infty$, $\mathbb{L}^{p}$ is a Banach space under
the norm $\parallel X\parallel_{p}:=(\mathbb{E}[\mid
X\mid^{p}])^{\frac{1}{p}}$.

Similarly, we can define space $\mathbb{L}^{p}(\Omega_{t})$.  And we denote by
$\mathbb{L}^{p}_{b}(\Omega_{t})$ the completion of $B_{b}(\Omega_{t})$ and $\mathbb{L}^{p}_{c}(\Omega_{t})$
the completion of $C_{b}(\Omega_{t})$ under norm $\parallel\cdot\parallel_{p}=\mathbb{E}[\vert \cdot\vert^{p}]$, $0\leq t\leq \infty$.

Denis, et al.\cite{Denis et al. (2011)} proved that $\mathbb{L}^{1}_{c}(\Omega_{t})=L^{1}_{G}(\Omega_{t})\subset \mathbb{L}^{1}_{b}(\Omega_{t})$.

$\mathbb{L}^{\infty}:=\mathcal{L}^{\infty}/\mathcal{N}$ is a Banach
space under the norm
$$
\parallel X\parallel_{\infty}:=\inf\{M\geq 0: \mid X\mid\leq M,
q.s.\}.
$$

\begin{lemm}\label{finite}
If for any $p>0$, $X\in \mathbb{L}^{p}$, then $c(\vert X\vert=\infty)=0$, that is $\vert X\vert<\infty, q.s.$
\end{lemm}
Proof. Otherwise, $c(\vert X\vert=\infty)>0$, then
\begin{equation*}
\mathbb{E}[\vert X\vert^{p}]\geq \mathbb{E}[\vert X\vert^{p}\mathbf{I}_{\{\vert X\vert=\infty\}}]=\infty,
\end{equation*}
which means $X\notin \mathbb{L}^{p}$.

Denote
\begin{equation*}
\begin{split}
\mathbb{L}_{S}(\Omega_{s,t}):=&\{\eta=\sum_{j=1}^{N}\mathbf{I}_{A^{j}}\eta^{j}, \text{ where }  \{A^{j}\}_{j=1}^{N} \text{ is an } \mathcal{F}_{s}-\text{partition of } \Omega,\\
 & \text{ and } \eta^{j}, j=1,\ldots,N \text{ are } \mathcal{F}^{s}_{t}-\text{measurable }\}.
\end{split}
\end{equation*}

\begin{defi} We define a mapping, $\mathbb{E}_{s}[\cdot]: \mathbb{L}_{S}(\Omega_{s,t})\rightarrow \mathbb{L}_{S}(\Omega_{s,s})$ which has the following properties,

(i) if $\xi\geq \eta$, then $\mathbb{E}_{s}[\xi]\geq \mathbb{E}_{s}[\eta]$.

(ii) $\mathbb{E}_{s}[\eta]=\eta$, if $\eta\in \mathbb{L}_{S}(\Omega_{s,s})$.

(iii) $\mathbb{E}_{s}[\xi]-\mathbb{E}_{s}[\eta]\leq \mathbb{E}_{s}[\xi-\eta]$.

(iv) $\mathbb{E}_{s}[\eta \xi]=\eta^{+}\mathbb{E}_{s}[\xi]+\eta^{-}\mathbb{E}_{s}[-\xi]$,if $\eta\in \mathbb{L}_{S}(\Omega_{s,s})$.

(v) $\mathbb{E}_{s}[\eta]=\mathbb{E}[\eta]$, if $\eta$ is independent from $\mathcal{F}_{s}$.

Then for any $\eta\in \mathbb{L}_{S}(\Omega_{s,t})$, we have
\begin{equation*}
\mathbf{I}_{A^{j}}\mathbb{E}_{s}[\eta]=\mathbb{E}_{s}[\mathbf{I}_{A^{j}}\sum_{j=1}^{N}\mathbf{I}_{A^{j}}\eta^{j}]=\mathbb{E}_{s}[\mathbf{I}_{A^{j}}\eta^{j}]=\mathbf{I}_{A^{j}}\mathbb{E}[\eta^{j}].
\end{equation*}
Summarizing over $j$, we have
\begin{equation*}
\mathbb{E}_{s}[\eta]=\sum_{j=1}^{N}\mathbf{I}_{A^{j}}\mathbb{E}[\eta^{j}].
\end{equation*}
So we can define the conditional G-expectation of $\eta$ as
\begin{equation*}
\mathbb{E}_{s}[\eta]:=\mathbb{E}[\eta\mid \mathcal{F}_{s}]=\sum_{j=1}^{N}\mathbf{I}_{A^{j}}\mathbb{E}[\eta^{j}].
\end{equation*}
\end{defi}

Denote
\begin{equation*}
\begin{split}
\mathbb{L}_{S}(\Omega_{s,r,t}):=&\{\eta=\sum_{j=1}^{N}\mathbf{I}_{A^{j}_{s}}\eta^{j}_{s,r}\eta^{j}_{r,t}, \text{ where }  \{A^{j}_{s}\}_{j=1}^{N} \text{ is an } \mathcal{F}_{s}-\text{partition of } \Omega,\\
 &\eta^{j}_{s,r}, j=1,\ldots,N  \text{ are } \mathcal{F}^{s}_{r}-\text{measurable }, \eta^{j}_{r,t}, j=1,\ldots,N \text{ are } \mathcal{F}^{r}_{t}-\text{measurable }\}.
\end{split}
\end{equation*}

\begin{prop}
$\mathbb{E}_{s}[\cdot]: \mathbb{L}_{S}(\Omega_{s,r,t})\rightarrow \mathbb{L}_{S}(\Omega_{s,s})$ is dynamically consistent.
\end{prop}

Proof.
For any random variable $\eta\in\mathbb{L}_{S}(\Omega_{s,r,t})$,
we have
\begin{equation*}
\begin{split}
\mathbb{E}_{s}[\eta]&=\sum_{j=1}^{N}\mathbf{I}_{A^{j}_{s}}\mathbb{E}[\eta^{j}_{s,r}\eta^{j}_{r,t}]=\sum_{j=1}^{N}\mathbf{I}_{A^{j}_{s}}\mathbb{E}[\mathbb{E}[x\eta^{j}_{r,t}]_{x=\eta^{j}_{s,r}}]\\
                      &=\sum_{j=1}^{N}\mathbf{I}_{A^{j}_{s}}\mathbb{E}[\eta^{j}_{s,r}]\mathbb{E}[\eta^{j}_{r,t}].
\end{split}
\end{equation*}
And
\begin{equation*}
\mathbf{I}_{A^{j}_{s}}\mathbb{E}_{r}[\eta]=\mathbb{E}_{r}[\mathbf{I}_{A^{j}_{s}}\eta^{j}_{s,r}\eta^{j}_{r,t}]=\mathbf{I}_{A^{j}_{s}}\eta^{j}_{s,r}\mathbb{E}[\eta^{j}_{r,t}],
\end{equation*}
Summarizing over $j$, we have
\begin{equation*}
\mathbb{E}_{r}[\eta]=\sum_{j=1}^{N}\mathbf{I}_{A^{j}_{s}}\eta^{j}_{s,r}\mathbb{E}[\eta^{j}_{r,t}].
\end{equation*}
Hence
\begin{equation*}
\mathbb{E}_{s}[\mathbb{E}_{r}[\eta]]=\mathbb{E}_{s}[\sum_{j=1}^{N}\mathbf{I}_{A^{j}_{s}}\eta^{j}_{s,r}\mathbb{E}[\eta^{j}_{r,t}]]=\sum_{j=1}^{N}\mathbf{I}_{A^{j}_{s}}\mathbb{E}[\eta^{j}_{s,r}]\mathbb{E}[\eta^{j}_{r,t}].
\end{equation*}
So $\mathbb{E}_{s}[\cdot]$ is dynamically consistent.

By the proof of lemma 43 of Denis, et al.\cite{Denis et al. (2011)}, ``the collection of processes $(\theta_{t})_{t\in[s,T]}$ with
$\{\theta_{t}=\sum I_{A_{j}}\theta^{j}_{t}:\{A_{j}\}_{j=1}^{N}$ is an $\mathcal{F}_{s}$ partition of $\Omega$, $\theta^{j}$ is $(\mathbb{F}^{s})$-adapted\} is dense in $\mathcal{A}^{\Theta}_{s,T}$''. So for any indicator function $\mathbf{I}_{A}\in \mathcal{F}_{t}$,
there exist sequence $\zeta^{i}=\sum_{j=1}^{N^{i}_{1}}\mathbf{I}_{A^{ij}_{1}}\eta^{ij}, i=1,\ldots$,
where for every $i$, $\{A^{ij}_{1}\}_{j=1}^{N^{i}_{1}}$ is an $\mathcal{F}_{s}$-partition of $\Omega$ and
$\eta^{ij}$ are  $\mathcal{F}^{s}_{t}$-measuable, such that $\zeta^{i}\rightarrow \mathbf{I}_{A}, i\rightarrow \infty$.

While for any $\mathcal{F}_{t}$-measurable random variable $\eta$, there exists simple function sequence $\xi^{i}=\sum_{j=1}^{N^{i}_{2}}\mathbf{I}_{A^{ij}}\eta^{ij}, i=1,2,\ldots$ where $\{A^{ij}\}_{j=1}^{N^{i}_{2}}$ is an $\mathcal{F}_{t}$-partition of $\Omega$, and $\eta^{ij}$ are constants, such that $\xi^{i}\rightarrow \eta, i\rightarrow \infty$.

So $\mathbb{L}_{S}(\Omega_{s,t})$ is dense in $L^{0}(\Omega_{t})$, and $\mathbb{L}_{S}(\Omega_{s,r,t})$ is dense in $\mathbb{L}_{S}(\Omega_{s,t})$ as well as in $L^{0}(\Omega_{t})$.

Hence for any $\eta\in \mathbb{L}^{1}(\Omega_{t})$, there exists $\eta^{i}\in \mathbb{L}_{S}(\Omega_{s,t})$ such that $\eta^{i}\rightarrow \eta, i\rightarrow \infty$. We define the conditional G-expectation of $\eta$ as
\begin{equation}\label{conditional G-expectation}
\mathbb{E}_{s}[\eta]=\lim_{i\rightarrow \infty}\mathbb{E}_{s}[\eta^{i}].
\end{equation}

The conditional G-expectation $\mathbb{E}_{s}[\cdot]: \mathbb{L}^{1}(\Omega_{t})\rightarrow \mathbb{L}^{1}(\Omega_{s})$ defined in (\ref{conditional G-expectation}) has the following properties.

\begin{prop} For each $X,Y\in \mathbb{L}^{1}(\Omega)$,

(i) if $X\geq Y$, then $\mathbb{E}_{s}[X]\geq \mathbb{E}_{s}[Y]$.

(ii) $\mathbb{E}_{s}[\eta]=\eta$, if $\eta\in \mathbb{L}^{1}(\Omega_{s})$.

(iii) $\mathbb{E}_{s}[X]-\mathbb{E}_{s}[Y]\leq \mathbb{E}_{s}[X-Y]$.

(iv) $\mathbb{E}_{s}[\eta X]=\eta^{+}\mathbb{E}_{s}[X]+\eta^{-}\mathbb{E}_{s}[-X]$ for each bounded $\eta\in \mathbb{L}^{1}(\Omega_{s})$.

(v) $\mathbb{E}_{s}[\mathbb{E}_{t}[X]]=\mathbb{E}_{t\wedge s}[X]$, in particular, $\mathbb{E}[\mathbb{E}_{s}[X]]=\mathbb{E}[X]$.
\end{prop}

\section{It\^o's integral in $\mathbb{M}^{2}(0,T)$}

For $T\in R^{+}$, a partition $\pi_{T}$ of $[0, T]$ is a finite ordered subset $\pi_{T}=\{t_{0},t_{1},\ldots,t_{N}\}$ such
that $0=t_{0}<t_{1}<\ldots<t_{N}=T$. Let $\mu(\pi_{T}):=\max\{\vert t_{i+1}-t_{i}\vert: i=0,1,\ldots, N-1\}$, and
use $\pi^{N}_{T}=\{t^{N}_{0},t^{N}_{1},\ldots,t^{N}_{N}\}$ to denote a sequence of partitions of $[0,T]$
such that $\lim_{N\rightarrow \infty}\mu(\pi^{N}_{T})=0$.

Let $p\geq 1$ be fixed. We consider the following type of simple processes: for a given
partition $\pi_{T}=\{t_{0},t_{1},\ldots,t_{N}\}$ of $[0,T]$ we set
\begin{equation*}
\eta_{t}(\omega)=\sum_{k=0}^{N-1}\xi_{k}(\omega)\mathbf{I}_{[t_{k},t_{k+1})}(t),
\end{equation*}
where $\xi_{k}\in \mathbb{L}^{p}(\Omega_{t_{k}})$, $k=0,1,2,\ldots,N-1$ are given. The collection of these processes is denoted by $\mathbb{M}^{p,0}(0,T)$.

\begin{defi} For each $p\geq 1$, we denote by $\mathbb{M}^{p}(0,T)$ the completion of $\mathbb{M}^{p,0}(0,T)$ under the norm
\begin{equation*}
\parallel\eta\parallel_{\mathbb{M}^{p}(0,T)}:=\left\{\mathbb{E}[\int_{0}^{T}\vert\eta_{t}\vert^{p}dt]\right\}^{\frac{1}{p}}.
\end{equation*}
\end{defi}
It is clear that $\mathbb{M}^{p}(0,T)\supset \mathbb{M}^{q}(0,T)$ for $1\leq p\leq q$.

\begin{defi} For an $\eta\in \mathbb{M}^{p,0}(0,T)$ with $\eta_{t}(\omega)=\sum_{k=0}^{N-1}\xi_{k}(\omega)\mathbf{I}_{[t_{k},t_{k+1})}(t)$, the related
Bochner integral is
\begin{equation*}
\int_{0}^{T}\eta_{t}(\omega)dt:=\sum_{k=0}^{N-1}\xi_{k}(\omega)(t_{k+1}-t_{k}).
\end{equation*}
\end{defi}
For each $\eta\in \mathbb{M}^{p,0}(0,T)$, set
\begin{equation*}
\tilde{\mathbb{E}}_{T}[\eta]:=\frac{1}{T}\mathbb{E}[\int_{0}^{T}\eta_{t}(\omega)dt].
\end{equation*}
$\tilde{\mathbb{E}}_{T}: \mathbb{M}^{p,0}(0,T)\rightarrow \mathbb{R}$ forms a sublinear expectation, so under the natural
norm $\parallel\eta\parallel_{\mathbb{M}^{p}(0,T)}$, the Bochner integral can be extended from $\mathbb{M}^{p,0}(0,T)$ to $\mathbb{M}^{p}(0,T)$

We now introduce the definition of It\^o's integral. For simplicity, we first introduce
It\^o's integral with respect to 1-dimensional G-Brownian motion.

Let $(B_{t})_{t\geq 0}$ be a 1-dimensional G-Brownian motion with $G(\alpha)=\frac{1}{2}(\bar{\sigma}^{2}\alpha^{+}-\underline{\sigma}^{2}\alpha^{-})$,
where $0\leq \underline{\sigma}\leq \bar{\sigma}<\infty$.

\begin{defi} For each $\eta\in \mathbb{M}^{2,0}(0,T)$ with $\eta_{t}(\omega)=\sum_{k=0}^{N-1}\xi_{k}(\omega)\mathbf{I}_{[t_{k},t_{k+1})}(t)$, define
\begin{equation*}
I(\eta)=\int_{0}^{T}\eta_{t}dB_{t}:=\sum_{k=0}^{N-1}\xi_{k}(B_{t_{k+1}}-B_{t_{k}}).
\end{equation*}
\end{defi}

\begin{lemm} The mapping $I: \mathbb{M}^{2,0}(0,T)\rightarrow \mathbb{L}^{2}(\Omega_{T})$ is a continuous linear mapping
and thus can be continuously extended to $I: \mathbb{M}^{2}(0,T)\rightarrow \mathbb{L}^{2}(\Omega_{T})$, and we have
\begin{equation}\label{Ito integral extended 1}
\mathbb{E}[\int_{0}^{T}\eta_{t}dB_{t}]=0,
\end{equation}
\begin{equation}\label{Ito integral extended 2}
\mathbb{E}[(\int_{0}^{T}\eta_{t}dB_{t})^{2}]\leq \bar{\sigma}^{2}\mathbb{E}[\int_{0}^{T}\eta_{t}^{2}dt].
\end{equation}
\end{lemm}
Proof. Notice that $B_{t_{i+1}}-B_{t_{i}}$ is independent of $\mathcal{F}_{t_{i}}$, so for $\xi_{i}\in \mathbb{L}^{1}(\Omega_{t_{i}})$, we have
\begin{equation}
\mathbb{E}[\xi_{i}(B_{t_{i+1}}-B_{t_{i}})]=\mathbb{E}[-\xi_{i}(B_{t_{i+1}}-B_{t_{i}})]=0,
\end{equation}
and
\begin{equation}
\mathbb{E}[\xi^{2}_{i}(B_{t_{i+1}}-B_{t_{i}})^{2}-\bar{\sigma}^{2}\xi^{2}_{i}(t_{i+1}-t_{i})]=\mathbb{E}[\mathbb{E}_{t}[\xi^{2}_{i}(B_{t_{i+1}}-B_{t_{i}})^{2}-\bar{\sigma}^{2}\xi^{2}_{i}(t_{i+1}-t_{i})]]=0.
\end{equation}
Hence we get (\ref{Ito integral extended 1}) and  (\ref{Ito integral extended 2}) by the same procedure as Peng \cite{Peng (2010)}.

\begin{defi} We define, for a fixed $\eta\in \mathbb{M}^{2}(0,T)$, the stochastic integral
\begin{equation*}
\int_{0}^{T}\eta_{t}dB_{t}:=I(\eta).
\end{equation*}
\end{defi}
It is clear (\ref{Ito integral extended 1}) and  (\ref{Ito integral extended 2})  still holds for $\eta\in \mathbb{M}^{2}(0,T)$.

The It\^o's integral has the following properties,
\begin{prop} Let $\xi,\eta\in \mathbb{M}^{2}(0,T)$, and let $0\leq s\leq r\leq t\leq T$. Then we have
\begin{equation*}
\begin{split}
&(i) \int_{s}^{t}\eta_{u}dB_{u}=\int_{s}^{r}\eta_{u}dB_{u}+\int_{r}^{t}\eta_{u}dB_{u}, q.s.,\\
&(ii) \int_{s}^{t}(\alpha\eta_{u}+\theta_{u})dB_{u}=\alpha\int_{s}^{t}\eta_{u}dB_{u}+\int_{s}^{t}\theta_{u}dB_{u}, \text{ if }\alpha \in \mathbb{L}^{2}(\Omega_{s}),\\
&(iii) \mathbb{E}[X+\int_{s}^{T}\eta_{u}dB_{u}\mid \Omega_{s}]=\mathbb{E}[X\mid \Omega_{s}], \text{ for } X\in \mathbb{L}^{1}(\Omega).
\end{split}
\end{equation*}
\end{prop}

For the multi-dimensional case. Let $G(\cdot): \mathbb{S}(d)\rightarrow \mathbb{R}$ be a
given monotonic and sublinear function and let $(B_{t})_{t\geq0}$ be a d-dimensional G-Brownian motion. For each fixed $a\in \mathbb{R}^{d}$, we use $B^{a}_{t}:=\langle a, B_{t}\rangle$. Then
$(B^{a}_{t})_{t\geq0}$ is a 1-dimensional $G_{a}$-Brownian motion with $G_{a}(\alpha)=\frac{1}{2}(\sigma_{aa^{T}}^{2}\alpha^{+}-\sigma_{-aa^{T}}^{2}\alpha^{-})$,
where $\sigma^{2}_{aa^{T}}=2G(aa^{T})$ and  $\sigma^{2}_{-aa^{T}}=-2G(-aa^{T})$. Similar to 1-dimensional case, we can define It\^o's integral by
\begin{equation*}
I(\eta):=\int_{0}^{T}\eta_{t}dB^{a}_{t}, \text{ for } \eta\in \mathbb{M}^{2}(0,T).
\end{equation*}
and have (\ref{Ito integral extended 1}),  (\ref{Ito integral extended 2}) and Proposition 3.6.

Let $(\langle B\rangle_{t})_{t\geq 0}$ be the quadratic variation process of 1-dimensional G-Brownian motion.

Define a mapping:
\begin{equation*}
Q_{0,T}(\eta)=\int_{0}^{T}\eta_{t}d\langle B\rangle_{t}:=\sum_{j=0}^{N-1}\xi_{j}(\langle B\rangle_{t_{j+1}}-\langle B\rangle_{t_{j}}): \mathbb{M}^{1,0}(0,T)\rightarrow \mathbb{L}^{1}(\Omega_{T}).
\end{equation*}
We have the following lemma.

\begin{lemm}
For each $\eta\in \mathbb{M}^{1,0}(0,T)$,
\begin{equation}\label{simple quadratic inequality}
\mathbb{E}[\vert Q_{0,T}(\eta)\vert]\leq \bar{\sigma}^{2}\mathbb{E}[\int_{0}^{T}\vert\eta_{t}\vert dt].
\end{equation}
Thus $Q_{0,T}: \mathbb{M}^{1,0}(0,T)\rightarrow \mathbb{L}^{1}(\Omega_{T})$ is a continuous linear mapping. Consequently, $Q_{0,T}$ can
be uniquely extended to $\mathbb{M}^{1}(0,T)$, and we have
\begin{equation}\label{quadratic inequality}
\mathbb{E}[\vert \int_{0}^{T}\eta_{t}d\langle B\rangle_{t}\vert]\leq \bar{\sigma}^{2}\mathbb{E}[\int_{0}^{T}\vert\eta_{t}\vert dt], \forall \eta\in \mathbb{M}^{1}(0,T).
\end{equation}
\end{lemm}
Proof. Notice that $\mathbb{E}[\vert \xi_{j}\vert(\langle B\rangle_{t_{j+1}}-\langle B\rangle_{t_{j}})-\bar{\sigma}^{2}\vert \xi_{j}\vert(t_{j+1}-t_{j})]=0$.
Then it is easy to check that (\ref{simple quadratic inequality}) as well as (\ref{quadratic inequality})holds.

\begin{prop}\label{integral equality prop} For any $\eta\in \mathbb{M}^{2}(0,T)$, we have
\begin{equation}\label{integral equality}
\mathbb{E}[(\int_{0}^{t}\eta_{s}dB_{s})^{2}]=\mathbb{E}[\int_{0}^{T}\eta_{s}^{2}d\langle B\rangle_{s}].
\end{equation}
\end{prop}
Proof. For $\eta\in \mathbb{M}^{2,0}(0,T)$, it is easy to check that (\ref{integral equality}) holds. We can continuously extend the
above equality to the case $\eta\in \mathbb{M}^{2}(0,T)$ and get  (\ref{integral equality}).

Similar to Li and Peng(2011), we can prove the following proposition.
\begin{prop}
For any $\xi\in \mathbb{M}^{1}(0,T)$, $\eta\in \mathbb{M}^{2}(0,T)$, and $0\leq t\leq T$, $\int_{0}^{t}\xi_{s}ds$, $\int_{0}^{t}\eta_{s}dB_{s}$,
and $\int_{0}^{t}\xi_{s}d\langle B\rangle_{s}$ are well defined processes which are continuous in $t$ quasi-surely.
\end{prop}

A stopping time $\tau$ with respect to filtration $(\mathcal{F}_{t})$ is a mapping $\tau: \Omega\rightarrow [0,T]$ such that for every $t$,
$\{\omega: \tau(\omega)\leq t\}\in \mathcal{F}_{t}$.

\begin{lemm}\label{stopping times in M^{p}}
For each stopping time $\tau$ and $\eta\in \mathbb{M}^{p}(0,T)$, we have $\mathbf{I}_{[0,\tau)}, \mathbf{I}_{[0,\tau)}\eta\in \mathbb{M}^{p}(0,T)$.
\end{lemm}
Proof.
For a given stopping time $\tau$, let
\begin{equation*}
\tau_{n}=\sum_{k=1}^{n}t^{n}_{k}\mathbf{I}_{[t^{n}_{k-1}<\tau\leq t^{n}_{k})}+T\mathbf{I}_{[\tau>T]}.
\end{equation*}
Then
\begin{equation*}
\begin{split}
\mathbf{I}_{[\tau_{n},T)}(t)=&\mathbf{I}_{[\sum_{l=1}^{n}t^{n}_{l}\mathbf{I}_{[t^{n}_{l-1}<\tau\leq t^{n}_{l})}+T\mathbf{I}_{[\tau>T]},T)}(t)\\
                                        =&\sum_{l=1}^{n}\mathbf{I}_{[t^{n}_{l},T)}(t)\mathbf{I}_{[t^{n}_{l-1}<\tau\leq t^{n}_{l})}\\
                                        =&\sum_{l=1}^{n}\sum_{k=l}^{n-1}\mathbf{I}_{[t^{n}_{k},t^{n}_{k+1})}(t)\mathbf{I}_{[t^{n}_{l-1}<\tau\leq t^{n}_{l})}\\
                                        =&\sum_{k=1}^{n-1}\left(\sum_{l=1}^{k}\mathbf{I}_{[t^{n}_{l-1}<\tau\leq t^{n}_{l})}\right)\mathbf{I}_{[t^{n}_{k},t^{n}_{k+1})}(t).
\end{split}
\end{equation*}
Since $\sum_{l=1}^{k}\mathbf{I}_{[t^{n}_{l-1}<\tau\leq t^{n}_{l})}\in \mathbb{L}^{p}(\Omega_{t^{n}_{k}})$,
we have $\mathbf{I}_{[\tau_{n},T)}\in \mathbb{M}^{p,0}(0,T)$.

For any $\eta\in \mathbb{M}^{p}(0,T)$, there exists a sequence of simple processes
$\eta^{i}_{t}=\sum_{k}\xi^{i}_{k}\mathbf{I}_{[t^{i}_{k},t^{i}_{k+1})}(t)$ with $\xi^{i}_{k}\in \mathbb{L}^{p}(\Omega_{t^{i}_{k}})$
such that $\eta^{i}\rightarrow\eta, i\rightarrow \infty$ in $\mathbb{M}^{p}(0,T)$. Obviously, $\mathbf{I}_{[\tau_{n},T)}\eta^{i}\in\mathbb{M}^{p,0}(0,T)$.

It is easy to check that $\mathbf{I}_{[\tau_{n},T)}\eta^{i}\rightarrow\mathbf{I}_{[\tau_{n},T)}\eta, i\rightarrow\infty$ in $\mathbb{M}^{p}(0,T)$, which means
that $\mathbf{I}_{[\tau_{n},T)}\eta \in \mathbb{M}^{p}(0,T)$.

Now we prove $\mathbf{I}_{[\tau_{n},T)}\eta\rightarrow\mathbf{I}_{[\tau,T)}\eta$ in $\mathbb{M}^{p}(0,T)$. We have
\begin{equation*}
\begin{split}
&\mathbb{E}\left[\int_{0}^{T}\vert \mathbf{I}_{[\tau_{n},T)}(t)\eta_{t}-\mathbf{I}_{[\tau,T)}(t)\eta_{t}\vert^{p} dt\right]\\
\leq & C\mathbb{E}\left[\int_{0}^{T}\vert \mathbf{I}_{[\tau_{n},T)}(t)\eta_{t}-\mathbf{I}_{[\tau_{n},T)}(t)\eta^{i}_{t}\vert^{p} dt\right]+C\mathbb{E}\left[\int_{0}^{T}\vert \mathbf{I}_{[\tau_{n},T)}(t)\eta^{i}_{t}-\mathbf{I}_{[\tau,T)}(t)\eta^{i}_{t}\vert^{p} dt\right]\\
&+C\mathbb{E}\left[\int_{0}^{T}\vert \mathbf{I}_{[\tau,T)}(t)\eta^{i}_{t}-\mathbf{I}_{[\tau,T)}(t)\eta_{t}\vert^{p} dt\right].
\end{split}
\end{equation*}
Since $\eta^{i}\rightarrow\eta, i\rightarrow \infty$ in $\mathbb{M}^{p}(0,T)$, for any $\epsilon>0$, there exists $I$ such that when $i\geq I$, $\mathbb{E}\left[\int_{0}^{T}\vert \eta^{i}_{t}-\eta_{t}\vert^{p} dt\right]<\frac{\epsilon}{3}$. Hence for some fixed $i\geq I$, we have
\begin{equation*}
\begin{split}
C\mathbb{E}\left[\int_{0}^{T}\vert \mathbf{I}_{[\tau_{n},T)}(t)\eta_{t}-\mathbf{I}_{[\tau_{n},T)}(t)\eta^{i}_{t}\vert^{p} dt\right]<\frac{\epsilon}{3},
\end{split}
\end{equation*}
\begin{equation*}
\begin{split}
C\mathbb{E}\left[\int_{0}^{T}\vert \mathbf{I}_{[\tau,T)}(t)\eta^{i}_{t}-\mathbf{I}_{[\tau,T)}(t)\eta_{t}\vert^{p} dt\right]<\frac{\epsilon}{3},
\end{split}
\end{equation*}
and
\begin{equation*}
\begin{split}
C\mathbb{E}\left[\int_{0}^{T}\vert \mathbf{I}_{[\tau_{n},T)}(t)\eta^{i}_{t}-\mathbf{I}_{[\tau,T)}(t)\eta^{i}_{t}\vert^{p} dt\right]\leq \sum_{k}\vert\xi^{i}_{k}\vert^{p}\mu(\pi^{n}_{T})<\frac{\epsilon}{3}, \text{ for some fixed } i\geq I \text{ and n large enough }.
\end{split}
\end{equation*}
So $\mathbf{I}_{[\tau_{n},T)}\eta\rightarrow \mathbf{I}_{[\tau,T)}\eta$ in $ \mathbb{M}^{p}(0,T)$, which means $\mathbf{I}_{[\tau,T)}\eta \in \mathbb{M}^{p}(0,T)$,
and consequently $\mathbf{I}_{[0,\tau_{n})}\eta\rightarrow \mathbf{I}_{[0,\tau)}\eta$ in $ \mathbb{M}^{p}(0,T)$.
As a special case, $\mathbf{I}_{[0,\tau)}\in \mathbb{M}^{p}(0,T)$.

By lemma \ref{stopping times in M^{p}}, the integral $\int_{0}^{t}\mathbf{I}_{[0,\tau]}(s)dB_{s}$ and $\int_{0}^{t}\mathbf{I}_{[0,\tau]}(s)\eta_{s}dB_{s}$ for $\eta\in \mathbb{M}^{2}(0,T)$ is well defined.

\begin{lemm}\label{stopping times integral}
For each stopping time $\tau$ and $\eta\in \mathbb{M}^{p}(0,T)$, we have
\begin{equation}\label{stopping times integral}
\int_{0}^{t\wedge\tau}\eta_{s}dB_{s}=\int_{0}^{t} \mathbf{I}_{[0,\tau)}\eta_{s}dB_{s}, q.s.
\end{equation}
\end{lemm}
Proof. Let
\begin{equation*}
\tau_{n}=\sum_{k=1}^{n}t^{n}_{k}\mathbf{I}_{[t^{n}_{k-1}<\tau\leq t^{n}_{k})}+T\mathbf{I}_{[\tau>T]}=\sum_{k=1}^{n}\mathbf{I}_{A^{k-1}_{n}}t^{n}_{k},
\end{equation*}
with $A^{k-1}_{n}=[t^{n}_{k-1}<\tau\leq t^{n}_{k})$, and $A^{n}_{n}=[\tau>T]$.

For any  $\eta\in \mathbb{M}^{p}(0,T)$, we have
\begin{equation*}
\int_{0}^{\tau_{n}}\eta_{s}dB_{s}=\int_{0}^{\sum_{k=1}^{n}\mathbf{I}_{A^{k-1}_{n}}t^{n}_{k}}\eta_{s}dB_{s}
=\sum_{k=1}^{n}\mathbf{I}_{A^{k-1}_{n}}\int_{0}^{t^{n}_{k}}\eta_{s}dB_{s}
=\int_{0}^{t}\sum_{k=1}^{n}\mathbf{I}_{[0,t^{n}_{k})}(s)\mathbf{I}_{A^{k-1}_{n}}\eta_{s}dB_{s}.
\end{equation*}
Thus we have
\begin{equation}\label{stopping times integral 1}
\int_{0}^{\tau_{n}}\eta_{s}dB_{s}=\int_{0}^{t}\mathbf{I}_{[0,\tau_{n})}(s)\eta_{s}dB_{s},
\end{equation}
and
\begin{equation*}
\int_{\tau_{n}}^{t\wedge\tau}\eta_{s}dB_{s}=\int_{0}^{t}\mathbf{I}_{[\tau_{n},t\wedge\tau)}(s)\eta_{s}dB_{s}\rightarrow 0,
\end{equation*}
by the continuity of $\int_{0}^{t}\eta_{s}dB_{s}$. Hence
\begin{equation}\label{stopping times integral 2}
\lim_{n\rightarrow \infty}\int_{0}^{\tau_{n}}\eta_{s}dB_{s}=\int_{0}^{t\wedge\tau}\eta_{s}dB_{s}-\lim_{n\rightarrow \infty}\int_{\tau_{n}}^{t\wedge\tau}\mathbf{I}_{[0,t\wedge\tau)}(s)\eta_{s}dB_{s}\rightarrow\int_{0}^{t\wedge\tau}\eta_{s}dB_{s}.
\end{equation}

By the proof of lemma \ref{stopping times integral}, $\mathbf{I}_{[0,\tau_{n})}\eta\rightarrow\mathbf{I}_{[0,\tau)}\eta, n\rightarrow\infty$,in $\mathbb{M}^{p}(0,T)$, so we have
\begin{equation*}
\int_{0}^{t}\mathbf{I}_{[0,\tau_{n})}(s)\eta_{s}dB_{s}\rightarrow\int_{0}^{t}\mathbf{I}_{[0,\tau)}(s)\eta_{s}dB_{s}, \text{ in } \mathbb{L}^{2}(\Omega_{t}).
\end{equation*}

By Denis, Hu and Peng(2011) proposition 17, there exists a subsequence $\int_{0}^{t}\mathbf{I}_{[0,\tau_{n_{k}}]}(s)\eta_{s}dB_{s}$ such that
\begin{equation}\label{stopping times integral 3}
\int_{0}^{t}\mathbf{I}_{[0,\tau_{n_{k}})}(s)\eta_{s}dB_{s}\rightarrow\int_{0}^{t}\mathbf{I}_{[0,\tau)}(s)\eta_{s}dB_{s}, q.s., \text{ as } k\rightarrow\infty,.
\end{equation}

From (\ref{stopping times integral 1}), (\ref{stopping times integral 2}), and (\ref{stopping times integral 3}), (\ref{stopping times integral}) holds.

\begin{lemm}\label{inequality of conditional G-expectation}
For $\eta\in \mathbb{M}^{2}(s,t)$ with $s<t$, we have
\begin{eqnarray}
\mathbb{E}_{s}[\vert\int_{s}^{t}\eta_{u}du\vert^{2}]\leq (t-s)\mathbb{E}_{s}[\int_{s}^{t}\eta^{2}_{u}du]\label{dt square},\\
\mathbb{E}_{s}[\vert\int_{s}^{t}\eta_{u}d\langle B\rangle_{u}\vert^{2}]\leq \bar{\sigma}^{4}(t-s)\mathbb{E}_{s}[\int_{s}^{t}\eta^{2}_{u}du]\label{quadratic variation square}.
\end{eqnarray}
\end{lemm}
Proof. For each fixed $\omega\in\Omega$, $\eta_{u}(\omega)$ is a measurable function on $[s,t]$. By lemma \ref{finite},
$\int_{s}^{t}\eta^{2}_{u}du<\infty, q.s.$. So for fixed $\omega\in\Omega$ such that $\int_{s}^{t}\eta^{2}_{u}(\omega)du<\infty$, we have
\begin{equation}
\vert\int_{s}^{t}\eta_{u}(\omega)du\vert^{2}\leq (\int_{s}^{t}\vert \eta_{u}(\omega)\vert du)^{2}\leq (t-s)\int_{s}^{t}\vert \eta_{u}(\omega)\vert^{2}du.
\end{equation}
Then we have
\begin{equation}
\mathbb{E}_{s}[\vert\int_{s}^{t}\eta_{u}dt\vert^{2}]\leq (t-s)\mathbb{E}_{s}[\int_{s}^{t}\vert \eta_{u}\vert^{2}dt].
\end{equation}

Now we prove (\ref{quadratic variation square}).
For $\eta^{n}_{u}=\sum_{i=0}^{n-1}\eta^{n}_{t_{i}}\mathbf{I}_{[t_{i},t_{i+1}]}(t)\in \mathbb{M}^{2,0}(0,T)$,
\begin{equation}
\begin{split}
&\mathbb{E}_{s}[\vert\int_{s}^{t}\eta^{n}_{u}d\langle B\rangle_{u}\vert^{2}]\\
=&\mathbb{E}_{s}[\vert\sum_{i=0}^{n-1}\eta^{n}_{t_{i}}(\langle B\rangle_{t_{i+1}}-\langle B\rangle_{t_{i}})\vert^{2}]\\
=&\mathbb{E}_{s}[\sum_{i,j=0}^{n-1}\vert \eta^{n}_{t_{i}}\eta^{n}_{t_{j}}\vert(\langle B\rangle_{t_{i+1}}-\langle B\rangle_{t_{i}})(\langle B\rangle_{t_{j+1}}-\langle B\rangle_{t_{j}})]\\
\leq&\bar{\sigma}^{4}\mathbb{E}_{s}[\sum_{i,j=0}^{n-1}\vert \eta^{n}_{t_{i}}\eta^{n}_{t_{j}}\vert(t_{i+1}-t_{i})(t_{j+1}-t_{j})]\\
=&\bar{\sigma}^{4}\mathbb{E}_{s}[(\int_{s}^{t}\eta^{n}_{u}dt)^{2}]\\
\leq&\bar{\sigma}^{4}(t-s)\mathbb{E}_{s}[\int_{s}^{t}\vert \eta^{n}_{u}\vert^{2}dt].
\end{split}
\end{equation}

Thus (\ref{quadratic variation square}) holds for $\eta^{n}\in \mathbb{M}^{2,0}(0,T)$. We can continuously extend the above equality to the case $\eta\in \mathbb{M}^{2}(0,T)$ and get (\ref{quadratic variation square}).

\section{It\^o's Formula}

\begin{lemm}\label{one step simple coefficient}
Let $\varphi\in C^{1,2}([0,T]\times \mathbb{R}^{n})$ with
$\partial_{t}\varphi, \partial_{x^{\nu}}\varphi, \partial_{x^{\mu}x^{\nu}}^{2}\varphi\in C_{b,L_{ip}}([0,T]\times \mathbb{R}^{n})$ for
$\mu,\nu=1,\ldots,n$. Let $s\in[0,T]$ be fixed and $X=(X^{1},\ldots,X^{n})^{T}$ be an n-dimensional process on $[s,T]$ of the form
\begin{equation*}
X^{\nu}_{t}=X^{\nu}_{s}+\alpha^{\nu}(t-s)+\eta^{\nu ij}(\langle B^{i}, B^{j}\rangle_{t}-\langle B^{i}, B^{j}\rangle_{s})+\beta^{\nu j}(B^{j}_{t}-B^{j}_{s}),
\end{equation*}
where for $\nu=1,\ldots, n$, $i,j=1,\ldots,d$, $\alpha^{\nu}, \eta^{\nu ij}\in \mathbb{L}^{4}(\Omega_{s})$, $\beta^{\nu j}\in \mathbb{L}^{8}(\Omega_{s})$ and $X_{s}=(X^{1}_{s},\ldots,X^{n}_{s})^{T}$
is a given random vector in $\mathbb{L}^{2}(\Omega_{s})$. Then for each $t\geq s$, we have, in $\mathbb{L}^{2}(\Omega_{t})$,
\begin{equation}\label{Ito C b Lip one step simple}
\begin{split}
\varphi(t,X_{t})-\varphi(s,X_{s})=&\int_{s}^{t}[\partial_{t}\varphi(u,X_{u})+\partial_{x^{\nu}}\varphi(u,X_{u})\alpha^{\nu}]du+\int_{s}^{t}\partial_{x^{\nu}}\varphi(u,X_{u})\beta^{\nu j} dB^{j}_{u}\\
                             +&\int_{s}^{t}[\partial_{x^{\nu}}\varphi(u,X_{u})\eta^{\nu ij}+\frac{1}{2}\partial^{2}_{x^{\mu}x^{\nu}}\varphi(u,X_{u})\beta^{\mu i}\beta^{\nu j}]d\langle B^{i}, B^{j}\rangle_{u}.
\end{split}
\end{equation}
Here we use the above repeated indices $\mu, \nu, i$ and $j$ imply the summation.
\end{lemm}
Proof. For each positive integer $N$, we set $\delta_{N}=(t-s)/N$ and take the partition
$$
\pi^{N}_{[s,t]}=\{t^{N}_{0},t^{N}_{1},\ldots,t^{N}_{N}\}=\{s,s+\delta_{N},\ldots,s+N\delta_{N}=t\}.
$$
We have
\begin{equation}\label{extended}
\begin{split}
\varphi(t,X_{t})-\varphi(s,X_{s})=&\sum_{k=0}^{N-1}[\varphi(t^{N}_{k+1},X_{t^{N}_{k+1}})-\varphi(t^{N}_{k},X_{t^{N}_{k}})]\\
                       =&\sum_{k=0}^{N-1}\{\partial_{t}\varphi(t^{N}_{k},X_{t^{N}_{k}})(t^{N}_{k+1}-t^{N}_{k})+\partial_{x^{\nu}}\varphi(t^{N}_{k},X_{t^{N}_{k}})(X^{\nu}_{t^{N}_{k+1}}-X^{\nu}_{t^{N}_{k}})\\
                        &+\frac{1}{2}[\partial^{2}_{x^{\mu}x^{\nu}}\varphi(t^{N}_{k},X_{t^{N}_{k}})(X^{\mu}_{t^{N}_{k+1}}-X^{\mu}_{t^{N}_{k}})(X^{\nu}_{t^{N}_{k+1}}-X^{\nu}_{t^{N}_{k}})+\kappa^{N}_{k}]\}
\end{split}
\end{equation}
where
\begin{equation*}
\begin{split}
\kappa^{N}_{k}=&\partial^{2}_{tt}\varphi(t^{N}_{k}+\theta\delta_{N},X_{t^{N}_{k}}+\theta(X^{\nu}_{t^{N}_{k+1}}-X^{\nu}_{t^{N}_{k}}))(t^{N}_{k+1}-t^{N}_{k})^{2}\\
               &+2\partial^{2}_{tx^{\nu}}\varphi(t^{N}_{k}+\theta\delta_{N},X_{t^{N}_{k}}+\theta(X^{\nu}_{t^{N}_{k+1}}-X^{\nu}_{t^{N}_{k}}))(t^{N}_{k+1}-t^{N}_{k})(X^{\nu}_{t^{N}_{k+1}}-X^{\nu}_{t^{N}_{k}})\\
               &+[\partial^{2}_{x^{\mu}x^{\nu}}\varphi(t^{N}_{k}+\theta\delta_{N},X_{t^{N}_{k}}+\theta(X^{\nu}_{t^{N}_{k+1}}-X^{\nu}_{t^{N}_{k}}))-\partial_{x^{\mu}x^{\nu}}\varphi(t^{N}_{k},X_{t^{N}_{k}})](X^{\nu}_{t^{N}_{k+1}}-X^{\nu}_{t^{N}_{k}})^{2}
\end{split}
\end{equation*}
with $\theta\in [0,1]$. We have, since $\partial_{x^{\mu}x^{\nu}}^{2}\varphi\in C_{b,L_{ip}}([0,T]\times \mathbb{R}^{n})$,
\begin{equation*}
\begin{split}
\mathbb{E}[\vert\sum_{k=0}^{N-1}\kappa^{N}_{k}\vert^{2}]\leq CN[\delta_{N}^{6}+\delta_{N}^{3}]\rightarrow 0,
\end{split}
\end{equation*}
where $C$ is a constant independent of $k$.

The rest terms in the summation of the right side of (\ref{extended}) are $\xi^{N}_{t}+\zeta^{N}_{t}$ with
\begin{equation*}
\begin{split}
\xi^{N}_{t}=&\sum_{k=0}^{N-1}\{\partial_{t}\varphi(t^{N}_{k},X_{t^{N}_{k}})(t^{N}_{k+1}-t^{N}_{k})\\
            &+\partial_{x^{\nu}}\varphi(t^{N}_{k},X_{t^{N}_{k}})[\alpha^{\nu}(t^{N}_{k+1}-t^{N}_{k})+\eta^{\nu ij}(\langle B^{i},B^{j}\rangle_{t^{N}_{k+1}}-\langle B^{i},B^{j}\rangle_{t^{N}_{k}})+\beta^{\nu j}(B^{j}_{t^{N}_{k+1}}-B^{j}_{t^{N}_{k}})]\\
            &+\frac{1}{2}\partial^{2}_{x^{\mu}x^{\nu}}\varphi(t^{N}_{k}, X_{t^{N}_{k}})\beta^{\mu i}\beta^{\nu j}(B^{i}_{t^{N}_{k+1}}-B^{i}_{t^{N}_{k}})(B^{j}_{t^{N}_{k+1}}-B^{j}_{t^{N}_{k}})\},
\end{split}
\end{equation*}
and
\begin{equation*}
\begin{split}
\zeta^{N}_{t}=&\frac{1}{2}\sum_{k=0}^{N-1}\{\partial^{2}_{x^{\mu}x^{\nu}}\varphi(t^{N}_{k}, X_{t^{N}_{k}})[\alpha^{\mu}(t^{N}_{k+1}-t^{N}_{k})+\eta^{\mu ij}(\langle B^{i},B^{j}\rangle_{t^{N}_{k+1}}-\langle B^{i},B^{j}\rangle_{t^{N}_{k}})]\\
            &\times [\alpha^{\nu}(t^{N}_{k+1}-t^{N}_{k})+\eta^{\nu ij}(\langle B^{i},B^{j}\rangle_{t^{N}_{k+1}}-\langle B^{i},B^{j}\rangle_{t^{N}_{k}})]\\
            &+2[\alpha^{\mu}(t^{N}_{k+1}-t^{N}_{k})+\eta^{\mu ij}(\langle B^{i},B^{j}\rangle_{t^{N}_{k+1}}-\langle B^{i},B^{j}\rangle_{t^{N}_{k}})]\beta^{\nu j}(B^{j}_{t^{N}_{k+1}}-B^{j}_{t^{N}_{k}})\}.
\end{split}
\end{equation*}

Now we prove $\xi^{N}_{t}$ converges to the right side of (\ref{Ito C b Lip one step simple}) and $\zeta^{N}_{t}$ converges to 0 in  $\mathbb{L}^{2}(\Omega_{t})$.

Firstly, we have the following estimates.
\begin{equation}\label{partial u varphi converge}
\begin{split}
&\mathbb{E}_{s}[\int_{s}^{t}\vert\partial_{u}\varphi(u,X_{u})-\sum_{k=0}^{N-1}\partial_{u}\varphi(t^{N}_{k},X_{t^{N}_{k}})\mathbf{I}_{[t^{N}_{k},t^{N}_{k+1})}(u)\vert^{2}du]\\
=&\mathbb{E}_{s}[\sum_{k=0}^{N-1}\int_{t^{N}_{k}}^{t^{N}_{k+1}}\vert\partial_{u}\varphi(u,X_{u})-\partial_{u}\varphi(t^{N}_{k},X_{t^{N}_{k}})\vert^{2}]du]\\
\leq&\sum_{k=0}^{N-1}\int_{t^{N}_{k}}^{t^{N}_{k+1}}\mathbb{E}_{s}[\vert\partial_{u}\varphi(u,X_{u})-\partial_{u}\varphi(t^{N}_{k},X_{t^{N}_{k}})\vert^{2}]du]\\
\leq&\sum_{k=0}^{N-1}\int_{t^{N}_{k}}^{t^{N}_{k+1}}C_{1}\mathbb{E}_{s}[\vert u-t^{N}_{k}\vert^{2}+\vert X_{u}-X_{t^{N}_{k}}\vert^{2}]du\\
\leq&\sum_{k=0}^{N-1}\int_{t^{N}_{k}}^{t^{N}_{k+1}}\mathbb{E}_{s}[(C_{1}+C_{2}(\alpha^{\nu})^{2})\vert u-t^{N}_{k}\vert^{2}+C_{2}(\eta^{\nu ij})^{2}\vert \langle B^{i},B^{j}\rangle_{u}-\langle B^{i},B^{j}\rangle_{t^{N}_{k}}\vert^{2}+C_{2}(\beta^{\nu j})^{2}\vert B^{j}_{u}-B^{j}_{t^{N}_{k}}\vert^{2}]du\\
\leq&(C_{1}+C_{2}(\alpha^{\nu})^{2})(t-s)\delta_{N}^{2}+C_{2}(\eta^{\nu ij})^{2}(t-s)\delta_{N}^{2}+C_{2}(\beta^{\nu j})^{2}(t-s)\delta_{N}.
\end{split}
\end{equation}
Similarly,
\begin{equation}\label{partial x varphi converge 1}
\begin{split}
&\mathbb{E}_{s}[\int_{s}^{t}\vert\partial_{x^{\mu}}\varphi(u,X_{u})-\sum_{k=0}^{N-1}\partial_{x^{\mu}}\varphi(t^{N}_{k},X_{t^{N}_{k}})\mathbf{I}_{[t^{N}_{k},t^{N}_{k+1})}(u)\vert^{2}du]\\
\leq&(C_{1}+C_{2}(\alpha^{\nu})^{2})(t-s)\delta_{N}^{2}+C_{2}(\eta^{\nu ij})^{2}(t-s)\delta_{N}^{2}+C_{2}(\beta^{\nu j})^{2}(t-s)\delta_{N},
\end{split}
\end{equation}
\begin{equation}\label{partial x varphi converge 2}
\begin{split}
&\mathbb{E}_{s}[\int_{s}^{t}\vert\partial_{x^{\mu}}\varphi(u,X_{u})-\sum_{k=0}^{N-1}\partial_{x^{\mu}}\varphi(t^{N}_{k},X_{t^{N}_{k}})\mathbf{I}_{[t^{N}_{k},t^{N}_{k+1})}(u)\vert^{2}d\langle B^{i},B^{j}\rangle_{u}]\\
\leq&(C_{1}+C_{2}(\alpha^{\nu})^{2})(t-s)\delta_{N}^{2}+C_{2}(\eta^{\nu ij})^{2}(t-s)\delta_{N}^{2}+C_{2}(\beta^{\nu j})^{2}(t-s)\delta_{N},
\end{split}
\end{equation}
and
\begin{equation}\label{partial xx varphi converge}
\begin{split}
&\mathbb{E}_{s}[\int_{s}^{t}\vert\partial_{x^{\mu}x^{\nu}}\varphi(u,X_{u})-\sum_{k=0}^{N-1}\partial_{x^{\mu}x^{\nu}}\varphi(t^{N}_{k},X_{t^{N}_{k}})\mathbf{I}_{[t^{N}_{k},t^{N}_{k+1})}(u)\vert^{2}d\langle B^{i},B^{j}\rangle_{u}]\\
\leq&(C_{1}+C_{2}(\alpha^{\nu})^{2})(t-s)\delta_{N}^{2}+C_{2}(\eta^{\nu ij})^{2}(t-s)\delta_{N}^{2}+C_{2}(\beta^{\nu j})^{2}(t-s)\delta_{N}.
\end{split}
\end{equation}

Of course, (\ref{partial u varphi converge})-(\ref{partial xx varphi converge}) implies
that $\sum_{k=0}^{N-1}\partial_{\cdot}\varphi(t^{N}_{k}, X_{t^{N}_{k}})\mathbf{I}_{[t^{N}_{k},t^{N}_{k+1})}(\cdot)$ converges
to $\partial_{\cdot}\varphi(\cdot, X_{\cdot})$, $\sum_{k=0}^{N-1}\partial_{x^{\mu}}\varphi(t^{N}_{k}, X_{t^{N}_{k}})\mathbf{I}_{[t^{N}_{k},t^{N}_{k+1})}(\cdot)$
converges to $\partial_{x^{\mu}}\varphi(\cdot, X_{\cdot})$,
and $\sum_{k=0}^{N-1}\partial^{2}_{x^{\mu}x^{\nu}}\varphi(t^{N}_{k}, X_{t^{N}_{k}})\mathbf{I}_{[t^{N}_{k},t^{N}_{k+1})}(\cdot)$ converges
to $\partial^{2}_{x^{\mu}x^{\nu}}\varphi(\cdot, X_{\cdot})$ in $\mathbb{M}^{2}(0,T)$.

Then
\begin{equation}\label{approximation 1}
\begin{split}
&\mathbb{E}[\vert\int_{s}^{t}[\partial_{u}\varphi(u,X_{u})+\partial_{x^{\nu}}\varphi(u,X_{u})\alpha^{\nu}]du+\partial_{x^{\nu}}\varphi(u,X_{u})\beta^{\nu j} dB^{j}_{u}\\
&+[\partial_{x^{\nu}}\varphi(u,X_{u})\eta^{n,\nu ij}+\frac{1}{2}\partial^{2}_{x^{\mu}x^{\nu}}\varphi(u,X_{u})\beta^{\nu j}\beta^{\nu j}]d\langle B^{i},B^{j}\rangle_{u}-\xi^{N}_{t}\vert^{2}]\\
\leq&5\{\mathbb{E}[\vert\int_{s}^{t}(\partial_{u}\varphi(u,X_{u})-\sum_{k=0}^{N-1}\partial_{u}\varphi(t^{N}_{k}, X_{t^{N}_{k}})\mathbf{I}_{[t^{N}_{k},t^{N}_{k+1})}(u)du\vert^{2}]\\
&+\mathbb{E}[\vert\int_{s}^{t}(\partial_{x^{\nu}}\varphi(u,X_{u})-\sum_{k=0}^{N-1}\partial_{x^{\nu}}\varphi(t^{N}_{k}, X_{t^{N}_{k}})\mathbf{I}_{[t^{N}_{k},t^{N}_{k+1})}(u))\alpha^{\nu} du\vert^{2}]\\
&+\mathbb{E}[\vert\int_{s}^{t}(\partial_{x^{\nu}}\varphi(u,X_{u})-\sum_{k=0}^{N-1}\partial_{x^{\nu}}\varphi(t^{N}_{k}, X_{t^{N}_{k}})\mathbf{I}_{[t^{N}_{k},t^{N}_{k+1})}(u))\beta^{\nu j} dB^{j}_{u}\vert^{2}]\\
&+\mathbb{E}[\vert\int_{s}^{t}(\partial_{x^{\nu}}\varphi(u,X_{u})-\sum_{k=0}^{N-1}\partial_{x^{\nu}}\varphi(t^{N}_{k}, X_{t^{N}_{k}})\mathbf{I}_{[t^{N}_{k},t^{N}_{k+1})}(u))\eta^{\nu ij} d\langle B^{i},B^{j}\rangle_{u}\vert^{2}]\\
&+\mathbb{E}[\vert\int_{s}^{t}\frac{1}{2}(\partial_{x^{\mu}x^{\nu}}\varphi(u,X_{u})-\sum_{k=0}^{N-1}\partial_{x^{\mu}x^{\nu}}\varphi(t^{N}_{k}, X_{t^{N}_{k}})\mathbf{I}_{[t^{N}_{k},t^{N}_{k+1})}(u))\beta^{\mu i}\beta^{\nu j}] d\langle B^{i},B^{j}\rangle_{u}\vert^{2}]\}\\
\leq&5(t-s)\{\mathbb{E}[\mathbb{E}_{s}[\int_{s}^{t}(\partial_{u}\varphi(u,X_{u})-\sum_{k=0}^{N-1}\partial_{u}\varphi(t^{N}_{k}, X_{t^{N}_{k}})\mathbf{I}_{[t^{N}_{k},t^{N}_{k+1})}(u))^{2}du]]\\
&+\mathbb{E}[(\alpha^{\nu})^{2}\mathbb{E}_{s}[\int_{s}^{t}(\partial_{x^{\nu}}\varphi(u,X_{u})-\sum_{k=0}^{N-1}\partial_{x^{\nu}}\varphi(t^{N}_{k}, X_{t^{N}_{k}})\mathbf{I}_{[t^{N}_{k},t^{N}_{k+1})}(u))^{2} du]]\\
&+K\mathbb{E}[(\beta^{\nu j})^{2}\mathbb{E}_{s}[\int_{s}^{t}(\partial_{x^{\nu}}\varphi(u,X_{u})-\sum_{k=0}^{N-1}\partial_{x^{\nu}}\varphi(t^{N}_{k}, X_{t^{N}_{k}})\mathbf{I}_{[t^{N}_{k},t^{N}_{k+1})}(u))^{2} du]]\\
&+K\mathbb{E}[(\eta^{\nu ij})^{2}\mathbb{E}_{s}[\int_{s}^{t}\vert\partial_{x^{\nu}}\varphi(u,X_{u})-\sum_{k=0}^{N-1}\partial_{x^{\nu}}\varphi(t^{N}_{k}, X_{t^{N}_{k}})\mathbf{I}_{[t^{N}_{k},t^{N}_{k+1})}(u)\vert^{2}du]]\\
&+\frac{1}{2}K\mathbb{E}[(\beta^{\mu j}\beta^{\nu j})^{2}\mathbb{E}_{s}[\int_{s}^{t}\vert\partial_{x^{\mu}x^{\nu}}\varphi(u,X_{u})-\sum_{k=0}^{N-1}\partial_{x^{\mu}x^{\nu}}\varphi(t^{N}_{k}, X_{t^{N}_{k}})\mathbf{I}_{[t^{N}_{k},t^{N}_{k+1})}(u)\vert^{2}du]]\}.
\end{split}
\end{equation}
By estimates (\ref{partial u varphi converge})-(\ref{partial xx varphi converge}), and $\alpha^{\nu}, \eta^{\nu ij}\in \mathbb{L}^{4}(\Omega_{s})$, $\beta^{\nu j}\in \mathbb{L}^{8}(\Omega_{s})$ we can prove the right side of (\ref{approximation 1}) converges to
0 as $N\rightarrow \infty$.

By the boundedness of $\partial^{2}_{x^{\mu}x^{\nu}}\varphi$ and
\begin{equation*}
\begin{split}
&\mathbb{E}_{t^{N}_{k}}[(\langle B^{i},B^{j}\rangle_{t^{N}_{k+1}}-\langle B^{i},B^{j}\rangle_{t^{N}_{k}})^{2}(B^{j}_{t^{N}_{k+1}}-B^{j}_{t^{N}_{k}})^{2}]\\
=&\mathbb{E}[(\langle B^{i},B^{j}\rangle_{t^{N}_{k+1}}-\langle B^{i},B^{j}\rangle_{t^{N}_{k}})^{2}(B^{j}_{t^{N}_{k+1}}-B^{j}_{t^{N}_{k}})^{2}]\\
\leq&\left\{\mathbb{E}[(\langle B^{i},B^{j}\rangle_{t^{N}_{k+1}}-\langle B^{i},B^{j}\rangle_{t^{N}_{k}})]^{4}\mathbb{E}[(B^{j}_{t^{N}_{k+1}}-B^{j}_{t^{N}_{k}})^{4}]\right\}^{\frac{1}{2}},
\end{split}
\end{equation*}

we have
\begin{equation*}
\begin{split}
\mathbb{E}[(\zeta^{N}_{t})^{2}]\leq&N C_{3}\sum_{k=0}^{N-1}\{\mathbb{E}[\vert\partial^{2}_{x^{\mu}x^{\nu}}\varphi(t^{N}_{k}, X_{t^{N}_{k}})\vert^{2}[(\alpha^{\mu})^{4}(t^{N}_{k+1}-t^{N}_{k})^{4}+(\eta^{\mu ij})^{4}(\langle B^{i},B^{j}\rangle_{t^{N}_{k+1}}-\langle B^{i},B^{j}\rangle_{t^{N}_{k}})^{4}\\
                                   &+2(\alpha^{\nu}\beta^{\nu j})^{2}(t^{N}_{k+1}-t^{N}_{k})^{2}(B^{j}_{t^{N}_{k+1}}-B^{j}_{t^{N}_{k}})^{2}\\
                                   &+2(\eta^{\nu ij}\beta^{\nu j})^{2}(\langle B^{i},B^{j}\rangle_{t^{N}_{k+1}}-\langle B^{i},B^{j}\rangle_{t^{N}_{k}})^{2}(B^{j}_{t^{N}_{k+1}}-B^{j}_{t^{N}_{k}})^{2}\}\\
                               \leq&N C_{4}\sum_{k=0}^{N-1}\{\mathbb{E}[(\alpha^{\mu})^{4}(t^{N}_{k+1}-t^{N}_{k})^{4}+(\eta^{\mu ij})^{4}\mathbb{E}_{t^{N}_{k}}[(\langle B^{i},B^{j}\rangle_{t^{N}_{k+1}}-\langle B^{i},B^{j}\rangle_{t^{N}_{k}})^{4}]\\
                                   &+2(\alpha^{\nu}\beta^{\nu j})^{2}(t^{N}_{k+1}-t^{N}_{k})^{2}\mathbb{E}_{t^{N}_{k}}[(B^{j}_{t^{N}_{k+1}}-B^{j}_{t^{N}_{k}})^{2}]\\
                                   &+2(\eta^{\nu ij}\beta^{\nu j})^{2}\mathbb{E}_{t^{N}_{k}}[(\langle B^{i},B^{j}\rangle_{t^{N}_{k+1}}-\langle B^{i},B^{j}\rangle_{t^{N}_{k}})^{2}(B^{j}_{t^{N}_{k+1}}-B^{j}_{t^{N}_{k}})^{2}]]\}\\
                               \leq&N C_{5}\sum_{k=0}^{N-1}\{\mathbb{E}[(\alpha^{\mu})^{4}+\vert\eta^{\mu ij}\vert^{4})\delta_{N}^{4}+2((\alpha^{\nu}\beta^{\nu j})^{2}+\vert\eta^{\nu ij}\beta^{\nu j}\vert^{2})\delta_{N}^{3}]\}\\
                                  =&C_{5}(N\delta_{N})^{2}\mathbb{E}[(\alpha^{\nu})^{4}+\vert\eta^{\nu ij}\vert^{4})\delta_{N}^{2}+2((\alpha^{\nu}\beta^{\nu j})^{2}+\vert\eta^{\nu ij}\beta^{\nu j}\vert^{2})\delta_{N}]\}\\
                          \rightarrow& 0, \text{ in } \mathbb{L}^{2}(\Omega_{t}).
\end{split}
\end{equation*}
We then have proved (\ref{Ito C b Lip}).

\begin{lemm}
Let $\varphi\in C^{1,2}([0,T]\times \mathbb{R}^{n})$ such that
$\partial_{t}\varphi, \partial_{x^{\mu}}\varphi, \partial_{xx}^{2}\varphi\in C_{b,L_{ip}}(\mathbb{R}^{n})$.
Let $X=(X^{1},\ldots,X^{n})$ and
\begin{equation*}
X_{t}=X_{s}+\int_{s}^{t}\alpha^{\nu}_{u}du+\int_{s}^{t}\eta^{\nu ij}_{u}d\langle B^{i},B^{j}\rangle_{u}+\int_{s}^{t}\beta^{\nu j}_{u}dB^{j}_{u},
\end{equation*}
where $\alpha^{\nu}, \eta^{\nu ij}\in \mathbb{M}^{4}(0,T)$ and $\beta^{\nu j}\in \mathbb{M}^{8}(0,T)$. Then for each $t\geq s$, we have, in $\mathbb{L}^{2}(\Omega_{t})$,
\begin{equation}\label{Ito C b Lip}
\begin{split}
\varphi(t,X_{t})-\varphi(s,X_{s})=&\int_{0}^{t}[\partial_{t}\varphi(u,X_{u})+\partial_{x^{\nu}}\varphi(u,X_{u})\alpha^{\nu}_{u}]du+\int_{0}^{t}\partial_{x^{\nu}}\varphi(u,X_{u})\beta^{\nu j}_{u}dB^{j}_{u}\\
                             +&\int_{0}^{t}[\partial_{x^{\nu}}\varphi(u,X_{u})\eta^{\nu ij}_{u}+\frac{1}{2}\partial^{2}_{x^{\mu}x^{\nu}}\varphi(u,X_{u})\beta^{\nu i}_{u}\beta^{\nu j}_{u}]d\langle B^{i},B^{j}\rangle_{u}.
\end{split}
\end{equation}
\end{lemm}
Proof.
For $\alpha^{\nu}, \eta^{\nu ij}\in \mathbb{M}^{4}(0,T)$ and $\beta^{\nu j}\in \mathbb{M}^{8}(0,T)$, there exist sequences of simple
processes $\alpha^{m,\nu} \overset{\mathbb{M}^{4}}{\longrightarrow}  \alpha^{\nu}, \eta^{m, \nu ij}\overset{\mathbb{M}^{4}}{\longrightarrow}  \eta^{\nu ij}, \beta^{m,\nu j}\overset{\mathbb{M}^{8}}{\longrightarrow}  \beta^{\nu j}$
as $m\rightarrow \infty$.
Let
\begin{equation*}
X^{m}_{t}=X_{s}+\int_{s}^{t}\alpha^{m, \nu}_{u}du+\int_{s}^{t}\eta^{m, \nu ij}_{u}d\langle B^{i},B^{j}\rangle_{u}+\int_{s}^{t}\beta^{m, \nu j}_{u}dB^{j}_{u}, m=1,2,\ldots.
\end{equation*}
From lemma \ref{one step simple coefficient}, we have
\begin{equation}\label{Ito C b Lip n}
\begin{split}
\varphi(t,X^{m}_{t})-\varphi(s,X^{m}_{s})=&\int_{0}^{t}[\partial_{t}\varphi(u,X^{m}_{u})+\partial_{x^{\nu}}\varphi(u,X^{m}_{u})\alpha^{m, \nu}_{u}]du+\int_{0}^{t}\partial_{x^{\nu}}\varphi(u,X^{m}_{u})\beta^{m, \nu j}_{u}dB^{j}_{u}\\
                             +&\int_{0}^{t}[\partial_{x^{\nu}}\varphi(u,X^{m}_{u})\eta^{m, \nu ij}_{u}+\frac{1}{2}\partial^{2}_{x^{\mu}x^{\nu}}\varphi(u,X^{m}_{u})\beta^{m, \mu i}_{u}\beta^{m, \nu j}_{u}d\langle B^{i},B^{j}\rangle_{u}.
\end{split}
\end{equation}
For any $s\leq t\leq T$,
\begin{equation*}
\mathbb{E}[\sup_{s\leq t\leq T}\vert X_{t}-X^{m}_{t}\vert^{4}]\leq C\mathbb{E}\int_{s}^{T}[(\alpha^{m,\nu ij}_{u}-\alpha^{\nu ij}_{u})^{4}+\vert\eta^{m, \nu ij}_{u}-\eta^{\nu j}_{u}\vert^{4}+\vert\beta^{m, \nu j}_{u}-\beta^{\nu ij}_{u}\vert^{4}]du\rightarrow 0,
\end{equation*}
since $\alpha^{m,\nu} \overset{\mathbb{M}^{4}}{\longrightarrow}  \alpha^{\nu}, \eta^{m, \nu ij}\overset{\mathbb{M}^{4}}{\longrightarrow}  \eta^{\nu ij}, \beta^{m,\nu j}\overset{\mathbb{M}^{8}}{\longrightarrow}  \beta^{\nu j}$. So
\begin{equation*}
\mathbb{E}\int_{s}^{T}\vert X_{t}-X^{m}_{t}\vert^{4}dt\leq C\int_{s}^{T}\mathbb{E}[\sup_{s\leq t\leq T}\vert X_{t}-X^{m}_{t}\vert^{4}]dt\rightarrow 0.
\end{equation*}

Then
\begin{equation*}
\begin{split}
&\mathbb{E}[\int_{0}^{T}\vert\partial_{x^{\nu}}\varphi(u,X^{m}_{u})\alpha^{m, \nu}_{u}-\partial_{x^{\nu}}\varphi(u,X_{u})\alpha^{\nu}_{u}\vert^{2}du]\\
\leq & 2\{\mathbb{E}[\int_{0}^{T}(\partial_{x^{\nu}}\varphi(u,X^{m}_{u})-\partial_{x^{\nu}}\varphi(u,X_{u}))^{2}\vert \alpha^{m, \nu}_{u}\vert^{2}du]+\mathbb{E}[\int_{0}^{T}(\partial_{x^{\nu}}\varphi(u,X_{u}))^{2}\vert \alpha^{m,\nu}_{u}-\alpha_{u}\vert^{2}du]\}\}\\
\leq & 2 C\{\mathbb{E}[\int_{0}^{T}\vert X^{m}_{u}-X_{u}\vert^{2}\vert \alpha^{m, \nu}_{u}\vert^{2}du+\mathbb{E}[\int_{0}^{T}\vert \alpha^{m,\nu}_{u}-\alpha^{\nu}_{u}\vert^{2}du]\}\\
\leq &2C\{\mathbb{E}[\int_{0}^{T}\vert X^{m}_{u}-X_{u}\vert^{4}du]\mathbb{E}[\int_{0}^{T}\vert \alpha^{m, \nu}_{u}\vert^{4}du]\}^{\frac{1}{2}}+C\mathbb{E}[\int_{0}^{T}\vert \alpha^{m,\nu}_{u}-\alpha^{\nu}_{u}\vert^{2}du]\\
\rightarrow& 0,
\end{split}
\end{equation*}
where $C$ is a constant.
Then we have proved $\partial_{x^{\nu}}\varphi(\cdot,X^{m}_{\cdot})\alpha^{m, \nu}_{\cdot}\rightarrow\partial_{x^{\nu}}\varphi(\cdot,X_{\cdot})\alpha^{\nu}_{\cdot}$ in $\mathbb{M}^{2}(0,T)$. Similarly, we can prove that in $\mathbb{M}^{2}(0,T)$,
\begin{equation*}
\begin{split}
&\partial_{\cdot}\varphi(\cdot,X^{m}_{\cdot})\rightarrow\partial_{\cdot}\varphi(\cdot,X_{\cdot}),\\
&\partial_{x^{\nu}}\varphi(\cdot,X^{m}_{\cdot})\eta^{m, \nu ij}_{\cdot}\rightarrow\partial_{x^{\nu}}\varphi(\cdot,X_{\cdot})\eta^{\nu ij}_{\cdot},\\
&\partial_{x^{\nu}}\varphi(\cdot,X^{m}_{\cdot})\beta^{m, \nu j}_{\cdot}\rightarrow\partial_{x^{\nu}}\varphi(\cdot,X_{\cdot})\beta^{\nu j}_{\cdot},\\
&\partial_{x^{\mu}x^{\nu}}\varphi(\cdot,X^{m}_{\cdot})\beta^{m, \mu i}_{\cdot}\beta^{m, \nu j}_{\cdot}\rightarrow\partial_{x^{\mu}x^{\nu}}\varphi(\cdot,X_{\cdot})\beta^{\mu i}_{\cdot}\beta^{\nu,j}_{\cdot}.
\end{split}
\end{equation*}
We then pass to limit as $m\rightarrow \infty$ in both sides of (\ref{Ito C b Lip n}) to get (\ref{Ito C b Lip}).


\begin{lemm}
Let $\varphi\in C^{1,2}([0,T]\times \mathbb{R}^{n})$ such that
$\varphi, \partial_{t}\varphi, \partial_{x^{\mu}}\varphi, \partial_{xx}^{2}\varphi$ are bounded and uniformly continuous on $[0,T]\times \mathbb{R}^{n}$.
Let $X=(X^{1},\ldots,X^{n})$ and
\begin{equation*}
X_{t}=X_{s}+\int_{s}^{t}\alpha^{\nu}_{u}du+\int_{s}^{t}\eta^{\nu ij}_{u}d\langle B^{i},B^{j}\rangle_{u}+\int_{s}^{t}\beta^{\nu j}_{u}dB^{j}_{u},
\end{equation*}
where $\alpha^{\nu}, \eta^{\nu ij}\in \mathbb{M}^{4}(0,T)$ and $\beta^{\nu j}\in \mathbb{M}^{8}(0,T)$. Then for each $t\geq s$, we have, in $\mathbb{L}^{2}(\Omega_{t})$,
\begin{equation}\label{Ito bounded and uniformly continuous}
\begin{split}
\varphi(t,X_{t})-\varphi(s,X_{s})=&\int_{0}^{t}[\partial_{t}\varphi(u,X_{u})+\partial_{x^{\nu}}\varphi(u,X_{u})\alpha^{\nu ij}_{u}]du+\int_{0}^{t}\partial_{x^{\nu}}\varphi(u,X_{u})\beta^{\nu j}_{u}dB^{j}_{u}\\
                             +&\int_{0}^{t}[\partial_{x^{\nu}}\varphi(u,X_{u})\eta^{\nu ij}_{u}+\frac{1}{2}\partial^{2}_{x^{\mu}x^{\nu}}\varphi(u,X_{u})\beta^{\mu i}_{u}\beta^{\nu j}_{u}]d\langle B^{i},B^{j}\rangle_{u}.
\end{split}
\end{equation}
\end{lemm}
Proof. We take $\{\varphi_{m}\}_{m=1}^{\infty}$ such that, for each $m$, $\varphi_{m}$ and all its first order and second order derivatives are in
$C_{b,L_{ip}}([0,T]\times \mathbb{R}^{n})$ and such that, as $m\rightarrow \infty$,
$\varphi_{m}, \partial_{t}\varphi_{m}, \partial_{x^{\nu}}\varphi_{m}$, and $\partial_{x^{\mu}x^{\nu}}^{2}\varphi_{m}$ converge respectively to
$\varphi, \partial_{t}\varphi, \partial_{x^{\nu}}\varphi$, and $\partial_{x^{\mu}x^{\nu}}^{2}\varphi$ uniformly on $[0,T]\times \mathbb{R}^{n}$.
We then use the It\^o's formula (\ref{Ito C b Lip}) to $\varphi_{m}(t,X_{t})$, i.e.,
\begin{equation*}
\begin{split}
\varphi_{m}(t,X_{t})-\varphi_{m}(s,X_{s})=&\int_{0}^{t}[\partial_{t}\varphi_{m}(u,X_{u})+\partial^{x^{\nu}}\varphi_{m}(u,X_{u})\alpha^{\nu ij}_{u}]du+\int_{0}^{t}\partial_{x^{\nu}}\varphi_{m}(u,X_{u})\beta^{\nu j}_{u}dB^{j}_{u}\\
                             +&\int_{0}^{t}[\partial_{x^{\nu}}\varphi_{m}(u,X_{u})\eta^{\nu ij}_{u}+\frac{1}{2}\partial^{2}_{x^{\mu}x^{\nu}}\varphi_{m}(u,X_{u})\beta^{\mu i}_{u}\beta^{\nu j}_{u}d\langle B^{i},B^{j}\rangle_{u}.
\end{split}
\end{equation*}
Passing to the limit as $m\rightarrow \infty$, we get (\ref{Ito bounded and uniformly continuous}).

\begin{theo}
Let $\varphi\in C^{1,2}([0,T]\times \mathbb{R}^{n})$ and $X=(X^{1},\ldots,X^{n})$ with
\begin{equation*}
X_{t}=X_{0}+\int_{0}^{t}\alpha^{\nu}_{u}du+\int_{0}^{t}\eta^{\nu ij}_{u}d\langle B^{i},B^{j}\rangle_{u}+\int_{0}^{t}\beta^{\nu j}_{u}dB^{j}_{u},
\end{equation*}
where $\alpha^{\nu}, \eta^{\nu ij}\in \mathbb{M}^{4}(0,T)$ and $\beta^{\nu j}\in \mathbb{M}^{8}(0,T)$. Then for each $t\geq s$, we have, in $\mathbb{L}^{2}(\Omega_{t})$,
\begin{equation}\label{Ito formula}
\begin{split}
\varphi(t,X_{t})-\varphi(s,X_{s})=&\int_{0}^{t}[\partial_{t}\varphi(u,X_{u})+\partial_{x^{\nu}}\varphi(u,X_{u})\alpha^{\nu ij}_{u}]du+\int_{0}^{t}\partial_{x^{\nu}}\varphi(u,X_{u})\beta^{\nu j}_{u}dB^{j}_{u}\\
                             +&\int_{0}^{t}[\partial_{x^{\nu}}\varphi(u,X_{u})\eta^{\nu ij}_{u}+\frac{1}{2}\partial^{2}_{x^{\mu}x^{\nu}}\varphi(u,X_{u})\beta^{\mu i}_{u}\beta^{\nu j}_{u}]d\langle B^{i},B^{j}\rangle_{u}.
\end{split}
\end{equation}
\end{theo}
Proof. For simplicity, we only prove the case of $n=d=1$. We set
$$
\gamma_{t}=\vert X_{t}\vert+\int_{0}^{t}(\vert\beta_{u}\vert^{8}+\vert\beta_{u}\vert^{4}+\vert\beta_{u}\vert^{4})du,
$$
and for $k=1,2,\ldots$, $\tau_{k}=\inf\{t\geq 0; \gamma_{t} \geq k\}$. Let $\varphi_{k}\in C^{1,2}([0,T]\times \mathbb{R}^{n})$,
such that $\partial_{t}\varphi_{k}, \partial_{x^{\nu}}\varphi_{k}, \partial_{x^{\mu}x^{\nu}}\varphi_{k}$ are bounded and uniformly continuous
and $\varphi_{k}=\varphi$, for $\vert x\vert\leq 2k$ and $t\in[0,T]$. By lemma \ref{stopping times in M^{p}},
$\mathbf{I}_{[0,\tau_{k}]}\alpha, \mathbf{I}_{[0,\tau_{k}]}\eta, \mathbf{I}_{[0,\tau_{k}]}\beta\in \mathbb{M}^{2}(0,T)$, so we have
\begin{equation*}
X_{t\wedge\tau_{k}}=X_{0}+\int_{0}^{t}\mathbf{I}_{[0,\tau_{k}]}\alpha_{u}du+\int_{0}^{t}\mathbf{I}_{[0,\tau_{k}]}\eta_{u}d\langle B^{i},B^{j}\rangle_{u}+\int_{0}^{t}\mathbf{I}_{[0,\tau_{k}]}\beta_{u}dB^{j}_{u},
\end{equation*}
Then we can apply the It\^o's formula (\ref{Ito bounded and uniformly continuous}) to $\varphi_{k}(t, X_{t\wedge\tau_{k}})$ to obtain
\begin{equation*}
\begin{split}
\varphi_{k}(t,X_{t\wedge\tau_{k}})-\varphi_{k}(s,X_{s})=&\int_{0}^{t}[\partial_{t}\varphi_{k}(u,X_{u})+\partial_{x}\varphi_{k}(u,X_{u})\alpha_{u}]\mathbf{I}_{[0,\tau_{k}]}du\\
                                                        &+\int_{0}^{t}\partial_{x}\varphi_{k}(u,X_{u})\beta_{u}\mathbf{I}_{[0,\tau_{k}]}dB^{j}_{u}\\
                                                        &+\int_{0}^{t}[\partial_{x}\varphi_{k}(u,X_{u})\eta_{u}+\frac{1}{2}\partial^{2}_{xx}\varphi_{k}(u,X_{u})\beta^{2}_{u}]\mathbf{I}_{[0,\tau_{k}]}d\langle B^{i},B^{j}\rangle_{u}.
\end{split}
\end{equation*}
Passing to the limit as $k\rightarrow \infty$, we get (\ref{Ito formula}).

\end{document}